\documentclass[12pt,a4paper,reqno]{amsart}
\usepackage{amsmath,amsthm}
\usepackage{amsfonts}
\usepackage{amssymb}
\allowdisplaybreaks

\newcommand{\be}{\begin{equation}}
\newcommand{\ef}{\end{equation}}
\chardef\bslash=`\\ % p. 424, TeXbook
\newtheorem*{thm*}{Theorem}

\theoremstyle{definition}

\newtheorem*{remark*}{Remarks}
\newtheorem*{defn*}{Definition}

\theoremstyle{remark}
\newcommand{\G}{\Gamma}
\newcommand{\wt}{\widetilde}
\newcommand{\wh}{\widehat}
\newcommand{\fc}{\frac}

\newcommand{\iy}{\infty}
 \renewcommand{\sectionmark}[1]{}

\newcommand{\Te} {Teichm\"{u}ller}

 \usepackage{amsfonts}
\newcommand{\field}[1]{\mathbb{#1}}
\newcommand{\g}{\gamma}

\newcommand{\dl}{\delta}
\newcommand{\D}{\field{D}}

\newcommand{\z}{\zeta}
\newcommand{\ov}{\overline}
\newcommand{\vp}{\varphi}
\newcommand{\hC}{\wh{\field{C}}}
\newcommand{\C}{\field{C}}

\newcommand{\B}{\mathbf{B}}
\newcommand{\T}{\mathbf{T}}
\newcommand{\Belt}{\operatorname{Belt}}

\newcommand{\Fib}{\operatorname{Fib}}
\newcommand{\Teich}{\operatorname{Teich}}

\renewcommand{\a} {\alpha}

\begin{document}

\title{A new method for evaluation of polynomial coefficients}
\author{Samuel L. Krushkal}

\begin{abstract} 
This paper introduces a new approach to estimating polynomial coefficients, an important problem in complex analysis.
This approach is intrinsically connected with features of univalent functions and of Teichm\"{u}ller spaces.
\end{abstract}

\date{\today\hskip4mm({CoefPolynom(4).tex})}

\maketitle

\bigskip

{\small {\textbf {2020 Mathematics Subject Classification:} Primary: 30C10, 30F60;
Secondary 12D05, 30C50, 30C75, 31A05}

\medskip

}

\bigskip

\markboth{Samuel L. Krushkal}{A new method for evaluation of coefficients of polynomials} \pagestyle{headings}

\bigskip\bigskip
\centerline{\bf 1. INTRODUCTION AND MAIN THEOREMS}

\bigskip\noindent
{\bf 1.1}.
Estimating polynomial coefficients has a long history and involves a variety of methods, but still remains an important problem due to many applications. This topic was intensively investigated also in geometric function theory; see, e.g., \cite{Ba}, \cite{B}, \cite{Dz}, \cite{Ga}, \cite{Go}),
\cite{MMR}, \cite{Va}, \cite{Wa}.

The aim of this paper is to present a new approach to estimating, which provides more general results. This approach relies on the deep intrinsic features of univalent functions and of Teichm\"{u}ller spaces.

Denote the collection of polynomials
$$
p_n(z) = c_0 + c_1 z + \dots + c_n z^n
$$
of degree $n \ge 1$ by $P_n$; it is a linear space isomorphic to $C^{n+1}$. Consider in $P_n$ the  homogeneous (circularly invariant) bounded closed convex domains $G \subset P_n$ containing the origin, which means that $G$ contains with any $p_n$ its pre and post rotations
 \be\label{1}
p_n(z) \mapsto p_{n,\a, \beta}(z) =  e^{i \beta} p_n(e^{i \a} z)
\end{equation}
with independent $\a, \beta \in [- \pi, \pi]$. The Minkowski functional
 \be\label{2}
l_G(p) = \inf \big\{r: \ r^{-1} p \in G, \ r > 0\big\};
\end{equation}
of any such domain determines on $C^{n+1}$ the norm $\| \cdot \|_X =  l_G(p)$, which is equivalent to the Euclidean norm (moreover, since $P_n$ is an $(n + 1)$-dimensional domain, all norm on $P_n$ are equivalent), and domain $G$ is in this norm the (closed) unit ball of the indicated $(n + 1)$-dimensional Banach space $X$.

Now we distinguish the polynomial
  \be\label{3}
p_n^0(z) = c_0^0 + c_1^0 z + \dots + c_n^0 z^n \in G
\end{equation}
with maximal first coefficient on $G$, i.e., with
$$
|c_1^0| = \max_G |c_1(p_n)|.
$$

Our first main result is the following general theorem, which is useful for applications and provides
simultaneously an {\bf algorithm} for explicit calculation of upper bounds for coefficients.

\bigskip\noindent
{\bf Theorem 1}. {\it For any rotationally invariant bounded convex domain $G \subset P_n$ and any polynomial $p_n \in G$, its coefficients $c_{m_j}$ (including $m_j = 1$) are sharply estimated by the non-zero coefficients
$c_{m_j}^0$ of the extremal polynomial (3) maximizing $|c_1|$ on $G$:
 \be\label{4}
|c_{m_j}| \le |c_{m_j}^0|.
\end{equation}

This extremal $p_n^0$ is determined up to pre and post rotations (1).}

\bigskip
This theorem implies that {\it estimating the coefficients of arbitrary polynomials from domain $G$ actually
is reduced to the
evaluation of the value of}
$$
\max |c_1| = \max |p_n^\prime(0)| \quad \text{on} \ \ G,
$$
i.e., to a standard maximization problem on a convex domain in $\C^{n+1}$.

\bigskip\noindent
{\bf 1.2}.
Theorem 1 is a consequence of the following result.
Every polynomial $p_n$ can be regarded as the {\bf Schwarzian derivative}
$$
S_w(z) = \left(\frac{w^{\prime\prime}(z)}{w^\prime(z)}\right)^\prime
- \frac{1}{2} \left(\frac{w^{\prime\prime}(z)}{w^\prime(z)}\right)^2
$$
of a locally univalent holomorphic function $w(z) = a_0 + a_1 z+ \dots$ on the unit disk $\D = \{|z| < 1\}$,  obtained as the solution of the differential equation
 \be\label{5}
(w^{\prime\prime}(z)/w^\prime(z))^\prime - (w^{\prime\prime}(z)/w^\prime(z))^2/2 = p_n(z), \quad p_n \in P_n,
\end{equation}
satisfying the prescribed initial conditions
$w(0) = a_0, \ \ w^\prime(0) = a_1, \ \ w^{\prime\prime}(0) = a_2$ (equivalently, as ratio of two linearly independent solutions $\eta_1, \eta_2$ of the linear differential equation
$2 \eta^{\prime\prime} + f(z) \ \eta = 0$
satisfying $\eta_1(0) = 0, \ \eta_1^\prime(0) = 1$ and $\eta_2(0) = 1, \ \eta_2^\prime(0) = 0$).

Recall that the Schwarzians $S_w$ of all locally univalent functions on $\D$ belong to the complex Banach space $\B$ of hyperbolically bounded holomorphic functions on $\D$ with norm
$$
\|f\|_\B = \sup_\D (1 - |z|^2)^2 |f(z)|.
$$
This space is dual to the space $A_1(\D)$ of integrable holomorphic functions on $\D$, and $A_1(\D) \subset \B$.

The class of all univalent functions $w(z) = z + a_2 z^2 + \dots$ with $w(0) = 0, \ w^\prime(0) = 1$ is denoted by $S$; it is one of the main canonical classes in geometric function theory.

There is an intrinsic connection between the dense subclass $S_Q$ of $S$ formed by $w \in S$ admitting quasiconformal extension across the unit circle $\mathbb C^1$ and hence onto to the complementary disk
$\D^* = \{z \in \hC = \C \cup \{\iy\}: \ |z| > 1\}$. Namely, the Schwarzians $S_w$ of $w \in S_Q$ fill in the space $\B$ a bounded domain (containing the origin), which serves as a canonical model of the {\bf universal Teichm\"{u}ller space $\T$}.

Noting that for any Moebius map $\g \in SL(2, \C)$, we have the equalities
$$
S_{w_1 \circ \g}(z) = (S_{w_1} \circ \g) \g^\prime(z)^2, \quad S_{\g \circ w}(z) = S_w(z),
$$
one can regard every $S_w(z)$ as a quadratic differential $\vp = S_w(z) dz^2$ on $\D$.

Define for a distinguished rotationally invariant bounded domain $G \subset P_n$ the quantity
 \be\label{6}
a_2(G) = \sup \big\{|a_2(w)|: \ S_w \in G \big\}.
\end{equation}
In view of compactness of $S$ in the topology of locally uniform convergence on $\D$, the supremum in
(5) is attained on some function $w_0 \in S$. Its Schwarzian also is a polynomial of degree $n$,
\be\label{7}
S_{w_0}(z) = c_{0,0} + c_{m_1,0} z^{m_1} + \dots + c_{m_j,0} z^{m_j} + \dots + c_{n,0} z^n;
\end{equation}
all presented here coefficients
$$
c_{m_j,0} \ne 0, \quad \text{with} \ \ 1 < m_1 \le m_2 \le \dots \le m_{n-1} \le n.
$$

Using the Schwarzians, we shall prove the following theorem giving Theorem 1.

\bigskip\noindent
{\bf Theorem 2}. {\it Given a rotationally invariant bounded convex domain $G \subset P_n$ which has the common points with the boundary of $\T$, take its extremal polynomial $p_n^0$ maximizing the first coefficient $c_1$ and the univalent solutions $w(z)$ of the equation (5) with $p_n \in G$.
Let its invariant $a_2(G)$ be attained on a function $w_0 \in S$, whose Schwarzian $S_{w_0}$ is presented by (5), with nonzero  coefficients $c_{m_j,0}$.
Then the corresponding coefficients $c_{m_j}$ of any polynomial $p_n \in G$ are sharply estimated
by
 \be\label{8}
|c_{m_j}| \le |c_{m_j}^0| = |c_{m_j,0}|;
\end{equation}
hence the coefficients $c_{m_j,0}$ and $c_{m_j}^0$ with $m_j > 1$ do not vanish simultaneously.

If $c_{1,0} \ne 0$, then the estimate (8) is valid also for} $m_j = 1$.

\bigskip
The proof of Theorem 2 given below provides in the case $c_{1,0} \ne 0$ the equality
$$
p_n^0 = S_{w_0}
$$
up to pre and post rotations (1).

Note also that the corresponding linear differential equations
$$
2 \eta^{\prime\prime}(z) + p_n(z) \eta(z) = 0
$$
closely relate to special functions, and one has to find the ratios of
some combinations of such functions.
For example, the solutions of the simplest equation of such type
$2 \eta^{\prime\prime} + c_1 z \ \eta = 0$
are represented by linear combinations of cylindrical functions (cf. \cite{Ka}).

\bigskip\noindent
{\bf 1.3}. We mention two {\bf examples} related to above theorems.

Consider the collection $G$ of polynomials $p_n \in P_n$ with $p_n(0) = 0$ and $|p_n(z)| < M$ on $\D$
with fixed $M < \iy$.

This $G$ is an $n$-dimensional convex domain. By Schwarz's lemma, all coefficients of such $p_n$ satisfy
$|c_m| \le M$, and $\max_G |c_m| = M$ is attained on $p_m(z) = M z^m$ for any $m \le n$.

\bigskip
Another example is given by the collection of even polynomials
$$
p_{2n} = c_0 + c_2 z^2 + \dots + c_{2n} z^{2n}.
$$
Such polynomials form a complex $(n + 1)$-dimensional space. On this space, $c_1 \equiv 0$, so the arguments
applied below in the proof of Theorem 2 do not work.

\bigskip\noindent
{\bf 1.4}.
There is an interesting class of holomorphic functions formed by nonvanishing (zero free) functions.  Such functions arise and play an essential role in many fields of complex analysis, in particular, as the derivatives of locally univalent functions and of covering maps, and have been investigated by many authors.

There were several deep conjectures on the coefficients of nonvanishing holomorphic functions and their generalization by Hummel-Scheinberg-Zalcman problem for nonvanishing functions in the general Banach spaces.
(see, e.g., \cite{HSZ}, \cite{Kr2}-\cite{Kr4} and the references cited there).

The indicated Hummel-Scheinberg-Zalcman problem is to find the bound
$$
C_n(X) = \sup_{\mathcal B(X)} |c_n|,
$$
where $\mathcal B(X)$ denotes the class of functions $f(z) = \sum\limits_0^\iy c_n z^n, \ |z| < 1$, which belong to a given Banach space $X$ and satisfy $\|f\|_X \le 1, \ f(z) \ne 0$.

Recently this problem has been solved for a wide class of Hilbert spaces in \cite{Kr5}.

\bigskip
The following theorem solves this problem for Banach spaces formed by polynomials.
Consider again the rotationally invariant convex domains $G \subseteq P_n$ and denote the subset of nonvanishing polynomials from $G$ by $G_{\mathcal Z}$.

\bigskip\noindent
{\bf Theorem 3}. {\it For any rotationally invariant bounded convex subdomain $G$ of $P_n$ and all polynomials $p_n \in G$, whose zero are located outside of the disk $\D$, their coefficients are sharply estimated  similar to (8) with bounds given by a nonvanishing in $\D$ polynomial $p_n^0$, on which the maximum of $|c_1|$ on $G_{\mathcal Z}$ is attained.
}

\bigskip\noindent
{\bf 1.5}.
The main idea of the proof of the above theorems is the same as applied in \cite{Kr1}-\cite{Kr3}. It  involves lifting the functionals $J(w) = c_m$ onto the Teichm\"{u}ller space $\T$ and its cover $\T_1$. But now the situation is somewhat different, which is caused by absence for polynomials of quasiconfromal deformations applied in the indicated papers.
Also we essentially simplify one of crucial steps in the proof, which concerns subharmonity.

\bigskip
An essential role in the proof of such results belongs to the following lemma giving a useful generalization of Koebe's one-quarter theorem.

\bigskip\noindent
{\bf Lemma 1}. \cite{Kr3} {\it Let $w(z)$ be a univalent solution of the equation (5) satisfying $w(0) = 0, \ w^\prime(0) = e^{i \theta}$ with the fixed $\theta \in [-\pi, \pi]$.
and let $w_0(z) = e^{i \theta} z + a_2^0 z^2 + \dots$ be one of the maximizing functions for $a_2$. Then the image domain $w(\D)$ covers entirely the disk
$D_{1/(2 |a_2(G)|)} = \{|w| < 1/(2 |a_2(G)|)\}$.
The radius value $1/(2 |a_2(G)|)$ is sharp for this collection of functions, and the circle $\{|w| = 1/(2 |a_2(G)|)$ contains the points not belonging to $w(\D)$ if and only if $|a_2| = |a_2(G)|$
(i.e., when $w$ is one of the maximizing functions for $a_2$ on the set $G$).

The inverted functions
$$
W(z) = 1/w(1/z) = e^{-i \theta} z - a_2^0 + b_1 z^{-1} + b_2 z^{-2} + \dots
$$
with $\z \in \D^*$ map the disk $\D^*$ onto a domain whose boundary is entirely contained in the
disk $\{|W + a_2(G)| \le |a_2(G)|\}$.
}

\bigskip\bigskip
\centerline{\bf 2. REMARKS ON TEICHM\"{U}LLER SPACES}

\bigskip
The proofs of the above theorems essentially rely on the deep geometric and
analytic features of Teichm\"{u}ller spaces. We briefly recall the needed results on these spaces;
for details see, e.g., \cite{Be}, \cite{GL}, \cite{Le}.

It is technically more convenient to deal with functions from $\Sigma_Q$.

\bigskip\noindent
{\bf 2.1}.
The {\bf universal Teichm\"{u}ller space} $\T = \Teich(\D)$ is the space of quasisymmetric homeomorphisms of the unit circle $\mathbb S^1$ factorized by M\"{o}bius maps;  all Teichm\"{u}ller spaces have their biholomorphic copies in $\T$.

The canonical complex Banach structure on $\T$ is defined by factorization of the ball of the Beltrami coefficients (or complex dilatations)
$$
\Belt(\D)_1 = \{\mu \in L_\iy(\C): \ \mu|\D^* = 0, \ \|\mu\| < 1\},
$$
letting $\mu_1, \mu_2 \in \Belt(\D)_1$ be equivalent if the corresponding  quasiconformal maps $w^{\mu_1}, w^{\mu_2}$ (solutions to the Beltrami equation $\partial_{\ov{z}} w = \mu \partial_z w$
with $\mu = \mu_1, \mu_2$) coincide on the unit circle $\mathbb S^1 = \partial \D^*$ (hence, on $\ov{\D^*}$). Such $\mu$ and the corresponding maps $w^\mu$ are called $\T$-{\it equivalent}.

We shall use the following lemma, which allows one to use some other normalizations of
quasiconformally extendable functions.

\bigskip\noindent
{\bf Lemma 2}. \cite{Kr2} {\it For any Beltrami coefficient $\mu \in \Belt(\D^*)_1$ and any $\theta_0 \in [0, 2 \pi]$, there exists a point $z_0 = e^{i \a}$ located on $\mathbb S^1$ so that
$|e^{i \theta_0} - e^{i \a}| < 1$ and such that for any $\theta$ satisfying
$|e^{i \theta} - e^{i \a}| < 1$ the equation
$\partial_{\ov z} w =  \mu(z) \partial_z w$
has a unique homeomorphic solution $w = w^\mu(z)$, which is holomorphic on the unit disk $\D$
and satisfies
$$
w(0) = 0, \quad w^\prime(0) = e^{i \theta}, \quad w(z_0) = z_0.
$$
Hence, $w^\mu(z)$ is conformal and does not have a pole in $\D$ \ (so
$w^\mu(z_{*}) = \iy$ at some point $z_{*}$ with $|z_{*}| \ge 1$).  }

\bigskip
In particular, this lemma allows one to define the Teichm\"{u}ller spaces using the quasiconformally extendible  univalent functions $w(z)$ in the unit disk $\D$ normalizing these functions by
$$
w(0) = 0, \quad w^\prime(0) = e^{i \theta}, \quad w(1) = 1.
$$
All such functions are holomorphic in the disk $\D$.

\bigskip\noindent
{\bf 2.2}.
The Teichm\"{u}ller space $\T_1 = \Teich(\D_{*})$ {\bf of the punctured disk} $\D_{*} = \D \setminus \{0\}$ is formed by classes $[\mu]_{\T_1}$ of $\T_1$-{\bf equivalent} Beltrami coefficients $\mu \in \Belt(\D)_1$ so that the corresponding quasiconformal automorphisms $w^\mu$ of the unit disk coincide on both boundary components (unit circle $\mathbb S^1$ and the puncture $z = 0$) and are homotopic on $\D \setminus \{0\}$.
This space can be endowed with a canonical complex structure of a complex Banach manifold
and embedded into $\T$ using uniformization of $\D_{*}$ by a cyclic parabolic Fuchsian
group acting discontinuously on $\D$ and $\D^*$. The functions $\mu \in L_\iy(\D)$ are lifted to
$\D$ as the Beltrami measurable $(-1, 1)$-forms  $\wt \mu d\ov{z}/dz$ in $\D$ with respect to
$\G$, i.e., via $(\wt \mu \circ \g) \ov{\g^\prime}/\g^\prime = \wt \mu, \ \g \in \G$,
forming the Banach space $L_\iy(\D, \G)$; we extend these $\wt \mu$ by zero to $\D^*$. Then $\T_1$ is canonically isomorphic to the subspace $\T(\G) = \T \cap \B(\G)$,
where $\B(\G)$ consists of elements $\vp \in \B$ satisfying $(\vp \circ \g) (\g^\prime)^2 = \vp$
in $\D^*$ for all $\g \in \G$.

Due to {\bf the Bers isomorphism theorem}, {\it the space $\T_1$ is biholomorphically isomorphic
to the Bers fiber space
$$
\Fib(\T) = \{(\phi_\T(\mu), z) \in \T \times \C: \ \mu \in \Belt(\D)_1, \ z \in w^\mu(\D)\}
$$
over the universal Teichm\"{u}ller space $\T$ with holomorphic projection $\pi(\psi, z) = \psi$} (see \cite{Be}).

This fiber space is a bounded hyperbolic domain in $\B_2 \times \C$ and represents the collection of domains $D_\mu = w^\mu(\D)$ as a holomorphic family over the space $\T$. For every $z \in \D$,  its
orbit $w^\mu(z)$ in $\T_1$ is a holomorphic curve over $\T$.

The indicated isomorphism between $\T_1$ and $\Fib(\T)$ is induced by the inclusion map \linebreak
$j: \ \D_{*} \hookrightarrow \D$ forgetting the puncture at the origin via
$$
\mu \mapsto (S_{w^{\mu_1}}, w^{\mu_1}(0)) \quad \text{with} \ \
\mu_1 = j_{*} \mu := (\mu \circ j_0) \ov{j_0^\prime}/j_0^\prime,
$$
where $j_0$ is the lift of $j$ to $\D$.

The Bers theorem is valid for Teichm\"{u}ller spaces $\T(X_0 \setminus \{x_0\})$ of all punctured hyperbolic Riemann surfaces $X_0 \setminus \{x_0\}$; we use only its special case.

\bigskip\noindent
{\bf 2.3}.
The spaces $\T$ and $\T_1$ can be weakly (in the topology generated by the spherical metric on $\hC$) approximate by finite dimensional Teichm\"{u}ller spaces $\T(0, n)$ of punctured spheres (Riemann surfaces of genus zero)
$$
X_{\mathbf z} = \hC \setminus \{0, 1, z_1 \dots, z_{n-3}, \iy\}
$$
defined by ordered $n$-tuples $\mathbf z = (0, 1, z_1, \dots, z_{n-3}, \iy), \ n > 4$ with distinct
$z_j \in \C \setminus \{0, 1\}$ (the details see, e.g., in \cite{Kr3}).

Another canonical model of $\T(0, n) $ is obtained again using the uniformization. This space
is biholomorphic to a bounded domain in the complex Euclidean space $\C^{n-3}$.

Note also that all Teichm\"{u}ller spaces are complete metric spaces with intrinsic Teichm\"{u}ller metric defined by quasiconformal maps. By the Royden-Gardiner theorem, this metric equals the hyperbolic  Kobayashi metric determined by the complex structure (see, e.g., \cite{EKK}, \cite{GL}, \cite{Ro}).

\bigskip\bigskip
\centerline{\bf 3. PROOF OF THEOREM 2}

\bigskip
Using the linearity of $c_m(p_n)$ on $P_n$, we distinguish the domain (ball in $X$)
$$
G_\rho = \rho G
$$
of radius $\rho > 0$ chosen so that $G_\rho \subset \T$ and the boundaries of $\T$ have the common points. Note that  by the Ahlfors-Weill theorem \cite{AW}, any $f \in \B$ with norm $\|f\|_\B  = k < 2$ is the Schwarzian derivative $S_w = f$ of a function $w$ which is univalent on the disk $\D$ and admits $k$-quasiconformal extension across the unit circle $\mathbb S^1$ to $\hC$ with Beltrami coefficient
canonically connected with $f$.

\bigskip
The proof of the theorem will be accomplished in three steps.
We first assume that the Schwarzian $S_{w_0}$ of univalent solution $w_0$ of the equation (3) with maximal $|a_2|$ has nonzero first coefficient of $c_1(S_{w_0})$ and show that in this case $S_{w_0}$ must be equal to the extremal polynomial $p_n^0$ maximizing $|c_1|$ on $G$.

\bigskip\noindent
{\bf Step 1: Lifting the coefficient functionals onto Teichm\"{u}ller spaces}.
Consider the classes $S_{z_0, \theta}$ of univalent functions $w(z)$ on the unit disk $\D$ with expansions
$$
w(z) = a_1 z^{-1} + a_2 z^2 + \dots
$$
with
$a_1 = e^{i \theta}, \ - \pi < \theta \le \pi$, having the fix point at $z_0$ on the boundary unit circle  $\mathbb S^1$ and
admitting quasiconformal extensions across $\mathbb S^1$ to the whole Riemann sphere $\hC$ and take their union
$$
\wh S = \bigcup_{z_0 \in \mathbb S^1, \theta \in (-\pi,\pi]} S_{z_0,\theta}.
$$
The Beltrami coefficients $\mu_f(z) = \partial_{\ov z} f/\partial_z f$ of these extensions run over the unit ball
$$\Belt(D^*)_1 = \{\mu \in L_\iy(\C): \ \mu(z)|D = 0, \ \ \|\mu\|_\iy  < 1\}.
$$

Now pass to the inverted functions $W_w(z) = 1/w(1/z)$ which form the corresponding classes $\Sigma_{z_0,\theta}$ of nonvanishing univalent functions on the disk $\D^*$ with expansions
$$
W(z) =  e^{- i \theta} z + b_0 + b_1 z^{-1} + b_2 z^{-2} + \dots, \quad  F(z_0) = z_0,
$$
and let
$\Sigma^0 = \bigcup_{z_0,\theta} \Sigma_{z_0,\theta}$.

The coefficients $a_n$ of $f(z)$ and the corresponding coefficients $b_j$ of $F_f(z)$ are related by
$$
b_0 + e^{2i \theta} a_2 = 0, \quad b_n + \sum \limits_{j=1}^{n}
\epsilon_{n,j}  b_{n-j} a_{j+1} + \epsilon_{n+2,0} a_{n+2} = 0,
\quad n = 1, 2, ... \ ,
$$
where $\epsilon_{n,j}$ are the entire powers of $e^{i \theta}$. This
successively implies the representations of $a_n$ by $b_j$ via
 \be\label{9}
a_n = (- 1)^{n-1} \epsilon_{n-1,0}  b_0^{n-1} - (- 1)^{n-1} (n - 2)
\epsilon_{1,n-3} b_1 b_0^{n-3} + \text{lower terms with respect to} \ b_0.
\end{equation}
This transforms any coefficient functional $J(w)$ on $S$ depending on a finite set of  distinguished
coefficients $a_{m_1}, \dots, a_{m_s}$  into a coefficient functional $\wt J(W^\mu)$
on $\Sigma^0$ depending on the corresponding coefficients $b_j$ and lifts the functionals  $\wt J(W)$ and $J(w)$ holomorphically onto the universal Teichm\"{u}ller space $\T$ modelled via bounded domain in the space $\B$ of Schwarzians $S_W^\mu$. Our functional $J(w) = c_n$ is expressed in terms of these  coefficients.

To lift $J$ onto the covering space $\T_1$, we again pass to functional $\wh J(\mu) = \wt J(W^\mu)$
lifting $J$ onto the ball $\Belt(\D)_1$ and apply the $\T_1$-equivalence, i.e., the quotient map
$$
\phi_{\T_1}: \ \Belt(\D)_1 \to \T_1, \quad \mu \to [\mu]_{\T_1}.
$$
Thereby the functional $\wt J(W^\mu)$ is pushed down to a bounded holomorphic functional $\mathcal J$
on the space $\T_1$ with the same range domain.

Regarding by Bers' isomorphism theorem the points of $\T_1$ as the pairs
$X_{W^\mu} = (S_{W^\mu}, F^\mu(0))$, where $\mu \in \Belt(\D)_1$ obey $\T_1$-equivalence, we come to a holomorphic functional
  \be\label{10}
\mathcal J(X_{W^\mu}) = \mathcal J(S_{W^\mu}, \ t), \quad t = W^\mu(0).
\end{equation}
on $\T_1 = \Fib(\T)$ and have to investigate its restriction to the image in $\Fib (\T)$ of the original submanifold $G_n$.

By Lemma 1, the boundary of any domain $W^\mu(\D^*)$ for any $W^\mu(z) \in \Sigma^0$ is located in the disk $\{|W - b_0| \le a_2(G) \}$, and the second coordinate $t$ runs over some subdomain $D_\theta$ in the disk $\D_{2 a_2(G)} = \{|t| < 2 a_2(G) \}$ containing the origin.

\bigskip\noindent
{\bf Step 2: Subharmonicity and the range domain of $W^\mu(0)$}.
Now we show that the evaluation of the functional (13) on $\T_1$ is based on subharmonicity of some related functions.

By Zhuravlev's theorem, the intersection of $\T$ (and of $\T_1$) with any complex line is a union of at most a countable set of simply connected planar domains (see \cite{KK}, Part 1, Ch. 5; \cite{Zh}).

Hence, taking for every point $\vp \in \G$ the supremum of  $|\mathcal J(S_{W^\mu}, t)|$ component-wise over the indicated intersection (followed by upper semicontinuous regularization), one obtains the corresponding subharmonic function $u_\vp(t), \ t \in W^\mu(0)$. This yields that also the function
 \be\label{11}
\sup_{\vp \in G} \ u_\vp(t)
\end{equation}
is subharmonic in some domain $D$ in the disk $\D_{2 |a_2(G)|}$ containing the origin.

Finally, we have to establish the range domain of $W^\mu(0)$ for $S_{W^\mu}$ running over $G$ and describe the boundary points of this domain.

The features of construction of the Bers fiber space $\Fib(\T)$ imply that image of $\chi(G-\rho)$ in the space $\Fib(\T)$ also is a connected submanifold covering $\chi(G_\rho)$.

We select a dense subsequence $\{\theta_1, \theta_2, \dots\} \subset [-\pi, \pi]$, getting the corresponding coefficient functionals on the classes $S_{z_0,\theta_j}$ and $\Sigma_{z_0,\theta_j}$, and consider the sequence of increasing products of the quotient spaces
 \be\label{12}
\mathcal T_m = \prod_{j=1}^m \ \wh \Sigma_{\theta_j}/\thicksim \
= \prod_{j=1}^m \{(S_{W_{\theta_j}}, W_{\theta_j}^{\mu_j}(0)) \} \ \simeq \T_1 \times \dots \times \T_1,
\end{equation}
where the equivalence relation $\thicksim$ again means $\T_1$-equivalence.
The Beltrami coefficients  $\mu_j \in \Belt(\D)_1$ are chosen here independently. For any $\T_1$, presented in the right-hand side of (15), the corresponding values of $W_{\theta_j}^{\mu_j}(0)$ run over some domain $D_{\a_j} \subset \C$, and the corresponding collection
$\beta = (\beta_1, \dots, \beta_s)$
of the Bers isomorphisms
$\beta_j: \ \{(S_{W_{\theta_j}}, W_{\theta_j}^{\mu_j}(0))\} \to \Fib(\T)$
determines a holomorphic surjection of the space $\mathcal T_m$ onto the product of $m$ spaces
$\Fib(\T)$. This provides the corresponding holomorphic maps
$$
\mathbf J = (\mathcal J(S_{W_{\theta_1}}), \dots, \mathcal J(S_{W_{\theta_m}}))
$$
whose (polydisk) norm satisfies
$$
\max_{G^m} \|\mathbf J\| = \max_j |J(S_{W_{\theta_1}})| = \max_G |c_{m_j}|.
$$

Now pick an extremal polynomial $p_n^{0,m} \in G$ for
$$
|\mathcal J(S_{W^\mu}, t)| = |c_{m_j}^0| = \max_G |c_{m_j}(p_n)|,
$$
noting that by assumption $c_{m_j}^0 \ne 0$. Consider the corresponding equation (5) for  this $p_n^{0,m}$
and take its univalent solution $w_{0,m}(z)$.

Applying to this $w_{0,m}$ the above construction of quotient spaces (12), one obtains, in view of the circular symmetry of the limit space for (12),  a maximal function
$$
u(t) = \max |c_{m_j}(t)|
$$
on the disk $\D_\a$ of some radius $\a \le 2 a_2(G)$; this function is circularly symmetric and positive on any circle $\{|t| = r\}$ and attains its maximal value on the boundary circle.

Since this construction involves the polynomials and univalent $w(z)$ with $a_2(2)$ arbitrarily close to
$a_{2 a_2(G)}$, we have that $\a = 2 a_2(G)$. So, for every $c_m$,
 \be\label{13}
\max_{\D_{2 a_2(G)}} \ u(t) =  \sup_{\chi(G)} \mathcal J(S_{W^\mu}, t) = \max_G |c_{m_j}|
= \max_{|t| = 2 a_2(G)} \ u(t),
\end{equation}
provided that $c_{m_j,0} \ne 0, \ m_j \ge 1$.

Take in (5) the third initial condition with $a_2 = a_2(G)$. Then the equalities (13) and the uniqueness of solution of the Cauchy problem for the equation (5), together
with the rotational invariance, provide that every extremal polynomial $p_n^{0,m}$ must coincide with $S_{w_0}$ and, hence, is determined up to pre and post rotations (1).

This holds also for $m_j =1$. Therefore, we have one (up to the indicated rotations) extremal polynomial $p_n^0$
for all $m_j$ given by (3).

This completes the proof under assumption that $c_1(S_{w_0}) \ne 0$.

\bigskip\noindent
{\bf Step 3: \ Omitting the assumption of nonwanishing $c_{1,0}$}.
The collection of $p_n \in G_n$ with $c_1 = 0$ is an analytic set in $G$.
Take a $\dl$-neighborhood $U_{\dl}$ of this set with small $\dl > 0$ and delete it from $G$. Then $G_n \setminus \ov{U_{\dl}}$ is an open set in $G$. One can apply to its components Lemma 1 and other above arguments, getting the corresponding maximal subharmonic functions
  \be\label{14}
u_\dl(t) = \sup_{S_{W^\mu} \in G \setminus U_\dl} \ |\mathcal J(S_{W^\mu}, t)|
\end{equation}
and the corresponding boundary univalent functions $w_\dl$ maximizing these functionals, for which one
has the estimates of coefficients of type (8).

Take the sequences $\{\dl_s\}$ decreasing to zero. Then the supremum in (14) is monotone increasing as
$s \to \iy$, while the sets $G \setminus U_{\dl_s}$ exhaust increasingly $G$. So the extremal values
of the moduli of $c_1$ and $c_{m_j}$ increase.

On the other hand, the compactness of the class $S$ with respect to locally uniform convergence on $\D$
provides that the extremal univalent $w_s$ maximizing $a_2$ and $c_{m_j}$ on $G \setminus U_{\dl_s}$
converge to $w_0$. Note that simultaneously the corresponding polynomials $p_{n,\dl}^0$ maximizing
$c_1$ on $G \setminus U_\dl$ converge to $p_n^0$ extremal on the limit domain $G$.

The extremal values of $c_{m_j}$ on $G \setminus U_{\dl_s}$ being attained on $S_{w_s}$  approach the corresponding non-zero extremal values of $c_{m_j}$ on $G$, and these values satisfy (8).
This completes the proof of Theorem 2.

\bigskip
Note that the values $|c_1(S_{w_s})|$ can decrease to zero (while the limit of $c_1(p_{n,\dl}^0)$ does not be zero).

\bigskip\noindent
{\bf Corollary from the proof}. The above arguments also imply that if the given subdomain $G \subseteq P_n$ contains with every its point $p_n(z) = c_0 + c_1 z + \dots$ the truncated polynomials
$c_1 z^s, \ s \le n$, then {\it the extremal polynomial $p_n^0$ maximizing $|c_1|$ must satisfy
$$
|c_{m_{j+1}}^0| \ge |c_{m_j}^0| \ge |c_1^0|,
$$
so the collection of non-zero coefficients $c_{m_j}^0$ has the decreasing moduli}.

\bigskip\bigskip
\centerline{\bf 4. SECOND PROOF OF THEOREM 2}

\bigskip
There is another, more complicated, proof of Theorem 2 in the lines of \cite{Kr1}, \cite{Kr2}.
It concerns the subharmonicity of function
 \be\label{15}
\wt u_\theta(t) = \sup_{W^\mu} |\mathcal J(S_{W^\mu}, t)|,
\end{equation}
taking the supremum over all $S_{W^\mu} \in \T$ admissible for any fixed given $t = W^\mu(0)$, which runs
over some domain $D$. This proof involves the weak approximation of the space $\T$ and simultaneously of
its cover $\T_1$ by finite dimensional Teichm\"{u}ller spaces of the punctured spheres in the topology of convergence on $\hC$ in the spherical metric (as was indicated in Section $\bf 2$).

We briefly outline this proof, because the same arguments will be applied in the proof of Theorem 3.

Fix $\theta \in [-\pi, \pi]$ and, using the maps $F^\mu \in \Sigma_{Q,\theta}$,
apply a weak approximation of the underlying space $\T$ (and simultaneously of the space $\T_1$) by finite dimensional Teichm\"{u}ller spaces of the punctured spheres in the topology of locally uniform convergence
on $\C$.

Take the set of points
$$
E = \{e^{\pi s i/2^n}, \ s = 0, 1, \dots, 2^{n+1} - 1; \ n = 1, 2,
\dots\}
$$
(which is dense on the unit circle) and consider the punctured spheres
$$
X_m = \hC \setminus \{e^{\pi s i/2^n}, \ s = 0, 1, \dots, 2^{n+1} -
1\}, \quad m = 2^{n+1},
$$
and their universal holomorphic covering maps $g_m: \ \D \to X_m$ normalized by
$g_m(0) = 0, \ g_m^\prime(0) > 0$.

The radial slits from the infinite point to all the points $e^{\pi s i/2^n}$ form a canonical dissection $L_m$ of $X_m$ and define the simply connected surface $X_m^\prime = X_m \setminus L_m$. Any
covering map $g_m$ determines a Fuchsian group $\G_m$ of covering transformations uniformizing $X_m^\prime$, which act discontinuosly in both disks $\D$ and $\D^*$.

Every such group $\Gamma_m$ has a canonical (open) fundamental polygon $P_m$ of $\G_m$ in $\D$
corresponding to the dissection $L_m$. It is a regular circular $2^{n+1}$-gon centered at the origin
of the disk and can be chosen to have a vertex at the point $z = 1$.
The restriction of $g_m$ to $P_m$ is univalent, and as $m \to \iy$, these polygons entirely increase and exhaust the disk $\D$.

Similarly, we take in the complementary disk $\D^*$ the mirror polygons $P_m^*$ and the covering maps $g_m^*(z) = 1/\ov{g_m(1/ \ov z)}$ which define the mirror surfaces $X_m^*$.

Now we approximate the maps $F^\mu \in \Sigma_{Q, \theta}$ by homeomorphisms $W^{\mu_m}$ having in
$\D = \{|z| < 1\}$ the Beltrami coefficients
$$
\mu_m = [g_m]_{*} \mu := (\mu \circ g_m) \ov{g_m^\prime}/g_m^\prime, \ \ n = 1, 2, \dots \ .
$$
Each $W^{\mu_m}$ is again $k$-quasiconformal (where $k = \|\mu\|_\iy$)  and compatible with
the group $\G_m$. As $m \to \iy$, the coefficients $\mu_m$ are convergent to $\mu$ almost everywhere on $\C$; thus, the maps $W^{\mu_m}$ are convergent to $F^\mu$ uniformly in the spherical
metric on $\hC$.

Note also that $\mu_m$ depend holomorphically on $\mu$ as elements of $L_\iy$; hence, $F^{\mu_m}(0)$
is a holomorphic function of $t = W^\mu(0)$.

As a result, one obtains that the Beltrami coefficients
$$
\mu_{h,m} := [g_m]_{*} \mu_h
$$
and the corresponding values  $W^{\mu_{h,m}}(0)$ are holomorphic functions of the variable
$t = F^\mu(0)$.

By Hartogs theorem, the function $\mathcal J(S_{F^{\mu_m}}, t)$ with $t =
F^{\mu_m}(0)$ is jointly holomorphic in $(S_{F^{\mu_m}}, t) \in \mathcal F(\T)$.

We now choose in $\T(0, m) \setminus \{\mathbf 0\}$ represented as a subdomain of
the space $\B(\G_m)$ a countable dense subset
$$
E^{(m)} = \{\vp_1, \vp_2, \dots, \vp_p, \dots\}.
$$
For any of its point $\vp_p$,  the corresponding extremal
Teich\"{u}ller disk $\D(\vp_p)$ joining this point with the origin of $\B(\G_m)$ does
not meet other points from this set
(this follows from the uniqueness of Teichm\"{u}ller extremal map). Recall also that
each disk $\D(\vp_p)$ is formed by the Schwarzians $S_{W^{\tau \mu_{p;m}}}$ with
$|\tau| < 1$ and
$$
\mu_{p;m}(z) = |\psi_{p;m}(z)|/\psi_{p;m}(z)
$$
with appropriate $\psi_{p;m} \in A_1(\D, \G_m), \ \|\psi_{p;m}\|_1 = 1$.

The restrictions of the functional $\mathcal J(S_{F^{\tau \mu_{p;m}}}, t)$ to
these disks are holomorphic functions of $(\tau, t)$; moreover, the above construction
provides that all these restrictions are
holomorphic in $t$ in some common domain $D_m \subset \D_{2|a_2(G)|}$
containing the point $t = 0$, provided that $|\tau| \le k < 1$. We use the maximal
common holomorphy domain; it is located in a disk $\{|t| < r_0\}, \ r_0 < 2 a_2(G)$.

Maximization over $\tau$ implies the logarithmically subharmonic functions
$$
U_{p;m}(t) = \sup_{|\tau| <1} |\mathcal J(S_{W^{\tau \mu_{p;m}}}, t)| \quad
(t = F^{\mu_{p;m}}(0), \ \ p = 1, 2, \dots)
$$
in the domain $D_m$. We consider the upper envelope of this sequence
$$
u_m(t) = \sup_p U_{p;m}(t)
$$
defined in a domain $D_m$ containing the origin, and take its
upper semicontinuous regularization
$$
u_m(t) = \limsup\limits_{t^\prime \to t} u_m(t^\prime),
$$
which does not increase $\max |\mathcal J|$
(by abuse of notation, we shall denote the regularizations by the same letter as the original functions).

Repeating this for all $m$, one obtains the sequences of monotone increasing functions $u_m(t)$ and of increasing domains $D_m$ exhausting a domain $D_\theta = \bigcup_m D_m$
such that each $u_m$ is subharmonic on $D_m$. The limit function of this sequence is also subharmonic on
$D_\theta$ and equals the function (15).

\bigskip\bigskip
\centerline{\bf 5. PROOF OF THEOREM 3}

\bigskip
Pass to polynomials $p_{n,r}(z) = \fc{1}{r} \ p_n(r z)$ with $r < 1$ close to  $1$ and consider their collection $G_{1/r}$.
Each $p_{n,r}$ has a neighborhood $U(p_{n,r}, \epsilon(r))$ in $G_{1/r}$ filled by polynomials which are zero free in the disk $\D$.
The union
$$
G_{\mathcal Z,r} = \bigcup_{p_n \in \mathcal Z} U(p_{n,r}, \epsilon(r))
$$
of such neighborhoods is a rotationally invariant domain. One can apply to it the arguments from the
second proof of Theorem 2 using  instead of $a_2(G)$ the maximal value of $|a_2(w)|$ on the image of the set $G_{\mathcal Z}$ in $S$ (cf. \cite{Kr2}). In a similar fashion, one obtains a subharmonic function of type (15) representing the functionals $c_{m_j}(p_n)$ on $G_{\mathcal Z}$.
The corresponding extremal polynomials $p_{n,r}^0$ do not vanish in $\D$).

This provides for $p_{n,r}$ the bounds (4) and (8) on $G_{\mathcal Z,r}$ depending on $r$.
Letting $r \to 1$, one obtains in the limit the desired estimate of type (4).

\bigskip
{\small\it{ \leftline{Department of Mathematics, Bar-Ilan University, 5290002 Ramat-Gan, Israel} \leftline{and Department of Mathematics, University of Virginia,  Charlottesville, VA 22904-4137, USA}}

\end{document}